\documentclass[reqno,11pt]{amsart}

\usepackage{amsmath,amsfonts,amssymb,amsthm,epsfig}
\usepackage{mathtools}
\usepackage[english]{babel}
\usepackage{tikz}
\usepackage[colorlinks,linkcolor=red,citecolor=red]{hyperref}
\usetikzlibrary{arrows}
\usepackage{thm-restate}
\usepackage{comment}

 \allowdisplaybreaks
\numberwithin{equation}{section}

\def\R{\mathbb R}
\def\S{\mathbb S}
\def\N{\mathbb N}

\def\bal{\begin{aligned}}
\def\eal{\end{aligned}}
\def\proofof#1{\begin{proof}[Proof of #1]}
\def\Chi#1{\hbox{{\large $\chi$}{\Large $_{_{#1}}$}}}

\def\XXint#1#2#3{{\setbox0=\hbox{$#1{#2#3}{\int}$} \vcenter{\vspace{-1pt}\hbox{$#2#3$}}\kern-.5\wd0}}

\pgfarrowsdeclare{<<<}{>>>}
{
\arrowsize=0.2pt
\advance\arrowsize by .5\pgflinewidth
\pgfarrowsleftextend{-4\arrowsize-.5\pgflinewidth}
\pgfarrowsrightextend{.5\pgflinewidth}
}
{
\arrowsize=0.2pt
\advance\arrowsize by .5\pgflinewidth
\pgfpathmoveto{\pgfpointorigin}
\pgfpathlineto{\pgfpoint{-1.2mm}{-.8mm}}
\pgfusepathqstroke
\pgfpathmoveto{\pgfpointorigin}
\pgfpathlineto{\pgfpoint{-1.2mm}{.8mm}}
\pgfusepathqstroke
}

\newcounter{mt}

\def\maintheoremdeclaration#1{\stepcounter{mt}\newcounter{#1}\setcounter{#1}{\arabic{mt}}}

\maintheoremdeclaration{main}

\newtheorem{theorem}{Theorem}[section]
\newtheorem{lemma}[theorem]{Lemma}
\newtheorem{prop}[theorem]{Proposition}

\theoremstyle{definition}
\newtheorem{defin}[theorem]{Definition}
\newtheorem{remark}[theorem]{Remark}
\newtheorem*{ack}{Acknowledgments}

\title{The double spherical cap rearrangement of planar sets}

\author{C. Gambicchia}
\address[]{Scuola Normale Superiore
	\newline\indent
	piazza dei Cavalieri 7,
	56127 Pisa, Italy}
\email{chiara.gambicchia@sns.it}

\begin{document}

\begin{abstract}
In the theory of shape optimization the rearrangements of sets are a key concept, because they allow us to keep some properties of the original set while improving other aspects.
This paper is devoted to the proof of an isoperimetric property of the double spherical cap rearrangement of planar sets.
In particular, we prove that, under the assumption of disconnection of non-trivial spherical slices, the rearranged set has a lower perimeter than the original one.
In the general case, the symmetrized set does not decrease the perimeter, but we show that the ``excess'' is bounded above by $2\mathcal{H}^1(\Gamma)$, where $\Gamma$ denotes the set of radii such that the spherical slice is a non-trivial arc of circle.  
Additionally, the higher-dimensional case is briefly discussed; in particular, an explicit counterexample is given, thus explaining why an analogous result cannot hold. 
The main reason for this is that, in dimension $N=3$ or higher, the union of two spherical caps of equal size does not minimize the $(N-2)$-dimensional measure of the boundary.
\end{abstract}

\subjclass[2020]{28A75}
\keywords{spherical rearrangements, isoperimetric inequality, rigidity, symmetrization}

\maketitle

\section{Introduction}
A crucial tool in the theory of shape optimization is given by the rearrangements of sets, which allow us to keep some properties of the original set while improving other aspects.
A quite complete guide about this topic is surely contained in \cite{K} and in the references therein.
In particular, when symmetrizing a set in some way, it is often useful to know that the rearrangement has the same volume of the original set, while decreasing the perimeter; 
indeed, it is not rare the case where, in solving some variational problem, one would like to reduce the searching field to sets that have some symmetry property.
An iconical example is the isoperimetric problem, where a key step in the proof by De Giorgi has been the proof of a rigidity property of Steiner's inequality (see, for example \cite{DG}).
Lately, a very productive research field is developping around quantitative isoperimetric inequalities (see for instance \cite{Fu1,CL,BCH,GP} and the references therein). 
In this setting, a very useful symmetrization is given by the spherical rearrangement (first introduced in \cite{P}), which allows one to decrease the perimeter, while keeping the volume of the intersection between the set and all balls centered at the origin.  
This last symmetrization is widely studied in \cite{CPS}, where the authors address in particular the question of rigidity.
However, in some situations the spherical cap rearrangement is not useful. 
Indeed, this symmetrization does not keep the barycenter and this could be a problem.
This is the case, for instance, of problems involving the barycentric asymmetry, that is, the volume of the symmetric difference between a set and the ball with the same volume and the same barycenter (see for instance \cite{Fu1,Fu2,BCH,GP}).
The {\it double spherical cap rearrangement} is a different symmetrization that can be helpful in such situations.
Roughly speaking, while in the spherical rearrangement of a set $E$ all the $(N-1)$-dimensional area of $E\cap \partial B(r)$ is moved towards a fixed direction, say $e_N$, in the {\it double spherical cap rearrangement} this area is split in two equal parts, one of which is moved towards the direction $e_N$, while the other one is moved towards the direction $-e_N$.
The advantage of this rearrangement is that not only it keeps the volume of $E$, but also its barycenter, as soon as the procedure is centered at the barycenter of the set in question;
the disadvantage, however, is that the perimeter is not decreased in general. 
In this paper we will present this symmetrization and study its properties, showing in particular that the perimeter decreases in dimension $2$ if all non-trivial slices are disconnected.
\begin{remark}
    This kind of symmetrization was first introduced in \cite{B}, where Bonnesen works with convex sets.
    In particular, he proves that, by centering this rearrangement in the right point, the symmetrized set has a lower perimeter than the original set.
    This was used in \cite{C}, where the author refers to this construction as the ``Bonnesen semicircular symmetrization'', summarizing very well its main properties.
    In a way, this paper generalizes what Bonnesen proved. 
    Indeed, the ``special point'', where the symmetrization is centered, is the center of the annulus of minimal witdh containing the boundary of the set, which is actually such that all the non-trivial spherical slices are disconnected.
\end{remark}

\subsection*{Notation and setting}
\begin{figure}[!ht]
    \begin{center}
        \begin{tikzpicture}[scale=1.2]
            \draw[dotted] (0,0) circle (1);
            \draw[black, thick] (0.7071,0.7071) arc[start angle=45, end angle=-45, radius=1];
            \draw[black, thick] (-0.7071,0.7071) arc[start angle=135, end angle=225, radius=1];
            \filldraw (0.7071,0.7071) circle (0.75pt);
            \filldraw (0.7071,-0.7071) circle (0.75pt); 
            \filldraw (-0.7071,-0.7071) circle (0.75pt);
            \filldraw (-0.7071,0.7071) circle (0.75pt); 
            \draw (0,0) -- (0.7071,0.7071);
            \draw[dashed] (-1.5,0) -- (1.5,0);
            \begin{scope}
                \clip (0,0) -- (0.7,0.7) -- (0.7,0) -- (0,0);
                \draw (0,0) circle (0.3);
            \end{scope} 
            \node at (0.5,0.2) {$\theta$};
            \node at (0.3,0.5) {$r$};       
        \end{tikzpicture}
    \end{center}
    \caption{The sets $D_\theta(r)$ and $S_\theta(r)$.}
    \label{fig D(r) e S(r)}
\end{figure}
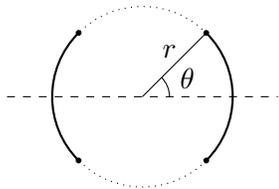

Let us start with some notation.
We will denote by $x=(x_1,x_2)$ a vector in $\mathbb{R}^2$ and by $\hat{x}=x/|x|$ the corresponding direction in $\S^1$.
Consider a set $E$ of finite perimeter and area. 
We denote by 
\[
\phi:\mathbb{R}^+\times \mathbb{S}^1\to \mathbb{R}^2
\]
the polar coordinates in the plane, \textit{i.e.} $\phi(r,\omega)=r\omega$.
We will sometimes prefer to refer to an angle $\theta\in [0,2\pi]$ rather than a direction $\omega\in \S^1$, so we set $\vartheta(\omega)$ as the (unique) angle such that 
\[
\omega=\big(\cos(\vartheta(\omega)), \sin(\vartheta(\omega))\big).
\]
For any positive radius $r$, we denote by $E_r=E\cap \partial B(r)$ the {\it spherical slice} of $E$, and we define the {\it circular distribution} of $E$ as 
\[
v_E(r)=\mathcal{H}^1(E_r).
\]
Notice that $v_E:(0, +\infty)\to [0,+\infty)$ is such that $v_E(r)\leq 2\pi r$ for any $r>0$.
For the sake of readability, however, we will drop the subscript $E$ in the notation.
Given a radius $r$ and an angle $\theta$ between $0$ and $\pi/2$, we denote by 
\[
D_{\theta}(r):= \left\{r\omega\in\partial B(r)\,|\,-\theta\leq \vartheta(\omega) \leq \theta\text{ or } \pi-\theta\leq \vartheta(\omega) \leq \pi+\theta\right\}
\]
the union of the two arcs of $\partial B(r)$, symmetrical with respect to the horizontal axis, of angle $2\theta$.
We then denote its ``$0$-dimensional boundary'' by $S_{\theta}(r)$, as in Figure \ref{fig D(r) e S(r)}.

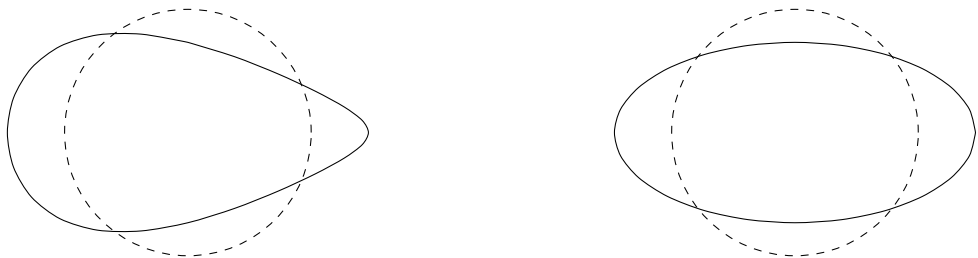
\begin{figure}[!h]
    \begin{center}
        \begin{tikzpicture}[scale=0.6]
        
            \draw[scale=2,domain=0:360,smooth,variable=\t] 
                plot ( {-2 * sin(\t)},{(1 + 0.5*sin(\t)) * cos(\t)});
            
            
            
            
            
            \draw[dashed] (0,0) circle (2.732);
            \end{tikzpicture}
            \hspace{3cm}
            \begin{tikzpicture}[scale=0.6]
                \draw[scale=2,domain=0:360,smooth,variable=\t] 
                    plot ({2 * cos(\t)}, {1 * sin(\t)});
                \draw[dashed] (0,0) circle (2.732);
            \end{tikzpicture}
    \end{center}
    \caption{A set and its double spherical cap rearrangement.}
    \label{fig double sph}
\end{figure}

\begin{remark}
    If $\theta=0$ or $\theta=\pi/2$, then $S_\theta(r)$ is empty; otherwise it is a union of four points, hence we have 
    \[
    \mathcal{H}^1(D_{\theta}(r))=4r\theta, \hspace{1.5cm} \mathcal{H}^0(S_{\theta}(r))=4\Chi{(0,\pi/2)}(\theta).
    \]
    Then, given $v$ as above, for any $r\in (0,+\infty)$ there exists a unique angle in $[0,\pi/2]$ such that 
    \[
    v(r)=\mathcal H^1(D_{\theta}(r)),
    \]
    and we denote it by $\theta_v(r)$. 
    In other words, we define 
    \begin{equation}\label{def theta_v}
        \theta_v(r):=\frac{\mathcal{H}^1(D_{\theta}(r))}{4r}=\frac{v(r)}{4r}.
    \end{equation} 
\end{remark}
We are now ready to define the double spherical cap rearrangement of a set, shown in Figure \ref{fig double sph}.
\begin{defin}
    Let $E$ be a set in $\mathbb{R}^2$ and let $v$ be its circular distribution.
    We define its {\em double spherical cap rearrangement} as 
    \[
    F_v:=\bigcup_{r>0}D_{\theta_v(r)}(r).
    \]
\end{defin} 

\begin{figure}[!ht]
    \begin{center}
        \begin{tikzpicture}
            \draw[dashed] (0,0) circle (1.2);
            \draw (0.9848,0.1736) arc[start angle=10, end angle=175, radius=1];
            \draw (0.9848,-0.1736) arc[start angle=350, end angle=185, radius=1];
            \draw (0.9848,0.1736) -- (2.9848,0.1736);
            \draw (2.9848,-0.1736) -- (0.9848,-0.1736);
            \draw (-0.9962,0.0871) -- (-4.9962,0.0871) -- (-4.9962,-0.0871) -- (-0.9962,-0.0871);
            \draw (2.9848,-0.1736) arc[start angle=-5, end angle=5, radius=1.9918];
        \begin{scope}[shift={(0,-2.6)}]
            \draw[dashed] (0,0) circle (1.2);
            \draw (0.9914,0.1305) arc[start angle=7.5, end angle=172.5, radius=1];
            \draw (0.9914,-0.1305) arc[start angle=352.5, end angle=187.5, radius=1];
            \draw (0.9914,0.1305) -- (2.9914, 0.1305) -- 
            (2.9914,0.0436) -- (4.9914, 0.0436) -- 
            (4.9914, -0.0436) -- (2.9914, -0.0436) --
            (2.9914,-0.1305) -- (0.9914,-0.1305);
            \draw (-0.9914,0.1305) -- (-2.9914, 0.1305) -- 
            (-2.9914,0.0436) -- (-4.9914, 0.0436) -- 
            (-4.9914, -0.0436) -- (-2.9914, -0.0436) --
            (-2.9914,-0.1305) -- (-0.9914,-0.1305);
        \end{scope}
        \end{tikzpicture}
    \end{center}
    \caption{The example showing that, in general, the double spherical cap symmetrization does not decrease the perimeter.}
    \label{fig counterex double sph}
\end{figure}
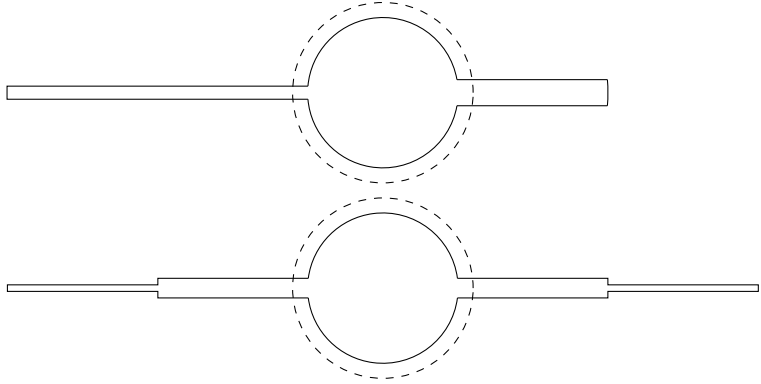

A few remarks are in order.
First of all, we notice that the double spherical rearrangement of any set is $2$-symmetrical by construction and consequently its barycenter is at the origin. 
Moreover, for every $r\geq 0$, this process does not affect the area of $E\cap B(r)$, and so this procedure does not affect the volume of the set.
It is then natural to wonder whether this procedure decreases the perimeter of the original set or not.
Unfortunately, the answer to the last question is that in general it does not.
Indeed, it is sufficient to consider a set made of a ball and two horizontal tentacles with different lengths but equal areas, so that the barycenter is at the center of the ball.
In this case, as one can see in Figure \ref{fig counterex double sph}, the symmetrization gives as a result a set of the same kind (a ball with two horizontal tentacles) but both tentacles of the rearranged set are as long as the longest tentacle of the original set, and so the perimeter is in fact increased.

In this example, the growth of the perimeter is due to the fact that for many radii the $r$-spherical slice $E_r$ consists in a single non trivial arc of circle, hence the ``split'' of such arcs in two pieces is not convenient in view of decreasing the perimeter. 
As soon as this does not happen, we can show that the double spherical cap symmetrization lowers the perimeter.
It is then useful to define the set of ``single-arched radii'' as follows
\begin{equation}\label{def Gamma_E}
\Gamma_E:=\{r\in (0,+\infty)\,|\, E_r \text{ is connected,  } 0<\mathcal{H}^1(E_r)<2\pi r\}.
\end{equation}
We can prove the following result.
\begin{theorem}\label{thm sph rearr}
    Let $E$ be a set in $\mathbb{R}^2$ with finite perimeter and finite volume and let $v:(0, +\infty)\to [0,+\infty)$ be its circular distribution.
    Then the following hold:
    \begin{itemize}
        \item $v$ is in $BV(0, +\infty)$;
        \item The double spherical cap rearrangment $F_v$ is a set of finite perimeter.
    \end{itemize}
    Moreover, defining $\Gamma_E$ as in \eqref{def Gamma_E}, one has
    \begin{equation}\label{strong ineq sph rearr}
        P\big(F_v,\phi(B\times \mathbb S^1)\big)\leq P\big(E,\phi(B\times \mathbb S^1)\big)+2\mathcal{H}^1(\Gamma_E\cap B)
    \end{equation}
    for any Borel set $B\subseteq (0,+\infty)$.
    In particular, if all non-trivial slices of $E$ are disconnected, then it holds 
    \begin{equation}\label{ineq sph rearr}
        P\big(F_v,\phi(B\times \mathbb S^1)\big)\leq P\big(E,\phi(B\times \mathbb S^1)\big)
    \end{equation}
    for any Borel set $B\subseteq (0,+\infty)$.
\end{theorem}
The proof of this statement makes use of standard techniques as, for example, those used in \cite{CPS}, where the authors consider the {\it spherical cap rearrangement} of sets.
The main difference is that the double spherical cap rearrangement does not, in general, decrease the perimeter of each slice of the set, hence we cannot use the isoperimetric inequality on the sphere.
However, the error is due to the radii in the set $\Gamma_E$ defined above; 
hence, by keeping track of these radii, we manage to obtain inequality \eqref{strong ineq sph rearr}.  

\subsection*{Organization of the paper}
In Section \ref{sec. prel.} the reader will find some preliminary results about BV-functions and sets of finite perimeter and some other results needed from the geometric measure theory.
Section \ref{sec. techn.} is devoted to the proof of a technical result regarding the map $\theta_v$ defined in \eqref{def theta_v} and a first estimate for the local perimeter of the rearrangement $F_v$.
Lastly, in Section \ref{sec. proof thm} we prove the main result of this work, giving a counterexample which explicits the reason why an analogous result cannot hold in the higher-dimensional case.

\begin{ack}
    The author would like to thank Giorgio Saracco, for pointing out references \cite{B, C}.
    The author is a member of the Gruppo Nazionale per l’Analisi Matematica, la Probabilità e le loro Applicazioni (GNAMPA) of the Istituto Nazionale di Alta Matematica (INdAM).
\end{ack}

\section{Preliminary notions}\label{sec. prel.}
\subsection*{Functions of bounded variation}
Let $f:\R^n\to \R$ be a Lebesgue-measurable function and let $\Omega$ be an open subset of $\R^n$.
We define the {\it total variation} of $f$ in $\Omega$ as 
\[
|Df|(\Omega):=\sup\left\{\int_\Omega f(x){\rm div}\varphi(x)dx\,\bigg\vert\,\varphi\in C^1_c(\Omega;\R^n)\,,\,|\varphi|\leq 1\right\}.
\]
We then say that $f$ belongs to $BV(\Omega)$ if $f\in L^1(\Omega)$ and it has bounded total variation in $\Omega$.
We shall write $f\in BV_{\rm loc}(\Omega)$ if, for any $\Omega'$ compactly contained in $\Omega$, one has $f\in BV(\Omega')$.
In particular, for $BV$-functions, the distributional derivative $Df$ is a vector-valued Radon measure, which can be decomposed as 
\[
Df=D^af+D^sf,
\]
where $D^af$ is absolutely continuous with respect to $\mathcal{H}^n$, with density denoted by $\nabla f$, while $D^sf$ is the singular part of the derivative.
Moreover, we can define the set of approximate-continuity points as those points for which there exists a $z\in\R$ such that
\[
\lim_{\rho\to 0^+}\frac{1}{|B_\rho(x)|}\int_{B_\rho(x)}|u(y)-z|dy=0.
\]  
Such $z$ is unique and is denoted by $\tilde u(x)$.
In particular, in the set $S_u$ of approximate-discontinuity points, we denote by $J_u$ the set of jump points, namely, those points for which there exist distinct $a,b\in\R$ and $\nu\in\S^{n-1}$ such that
\[
\int_{B_\rho(x)\cap H_\nu^+}|u(y)-a|dy=o(\rho^n),\hspace{1.5cm}\int_{B_\rho(x)\cap H_\nu^-}|u(y)-b|dy=o(\rho^n),
\]
where $H_\nu^+=\{\sigma\in\S^{n-1}| \sigma \cdot \nu>0\}$ and $H_\nu^-=\{\sigma\in\S^{n-1}| \sigma \cdot \nu<0\}$.
Such $a,b$ and $\nu$ are unique up to a change of sign of $\nu$ and a switch of $a$ and $b$; 
thus they are denoted, respectively, as $u^+(x), u^-(x)$ and $\nu_u(x)$. 
Specifically, $D^sf$ can in turn be decomposed as 
\[
D^sf=D^jf+D^cf,
\]
where the first term is the jump part of the derivative, while the second one is the Cantor part.
For brevity, we also write $\tilde{D}u=D^a u+D^c u$.
Among the fine properties of functions of bounded variations, we are interested in the following one (see \cite[Theorem 3.96 and Example 3.97]{AFP}).
\begin{prop}[Leibniz rule in BV]\label{Leibniz rule BV}
    For any pair $u_1,u_2$ of bounded functions in $BV(\Omega)$, we have that $u=u_1u_2\in BV(\Omega)$ with 
    \begin{equation}\label{eq. leibniz}
        \tilde{D}u=\tilde{u}_1\tilde{D}u_2+\tilde{u}_2\tilde{D}u_1 \hspace{1.5cm}D^j u=(u^+-u^-)\nu_u\mathcal{H}^{n-1}\llcorner J_u,
    \end{equation}
    where $J_u=J_{u_1}\cup J_{u_2}$, $\nu_u$ is consistently chosen to be equal to $\nu_{u_1}$ on $J_{u_1}\setminus S_{u_2}$, equal to $\nu_{u_2}$ on $J_{u_2}\setminus S_{u_1}$ and equal to $\nu_{u_1}=\nu_{u_2}$ on $J_{u_1}\cap J_{u_2}$, and where $u^+$ and $u^-$ are defined $\mathcal{H}^{n-1}$-$a.e.$ on $J_u$ as 
    \[
    u^+(x)=\begin{cases}
        \tilde{u}_2u^+_1& x\in J_{u_1}\setminus S_{u_2}\\
        \tilde{u}_1u^+_2& x\in J_{u_2}\setminus S_{u_1}\\
        u_2^+u^+_1& x\in J_{u_1}\cap J_{u_2},\\
    \end{cases}
    \hspace{1.5cm}
    u^-(x)=\begin{cases}
        \tilde{u}_2u^-_1& x\in J_{u_1}\setminus S_{u_2}\\
        \tilde{u}_1u^-_2& x\in J_{u_2}\setminus S_{u_1}\\
        u_2^-u^-_1& x\in J_{u_1}\cap J_{u_2}.\\
    \end{cases}
    \]
\end{prop}
\begin{remark}
    In particular, if $u_2$ is $C^1$ on $\Omega$, \eqref{eq. leibniz} can be zipped into 
    \[
        Du=\tilde{u}_1Du_2+u_2Du_1,
    \]
    and $J_u$ coincides with $J_{u_1}$.
\end{remark}

\subsection*{Sets of finite perimeter in $\R^2$}
\begin{figure}[!ht]
    \begin{center}
        \begin{tikzpicture}[scale=1.2]
                \begin{scope}
                    \clip (0,-2) rectangle (4,2);
                    \draw[thick] (0,0) .. controls (2,1) and (3,-2) .. (4.5,0) 
                                    .. controls (4.5,0) and (6,2) .. (8,0);
                    \draw[dotted] (4,0) circle (2);
                    \draw[->] (2,-0.1)-- (2.5, 0.8);
                    \draw[->] (2,-0.1) -- (2,0.8) node[left] {\small $\nu_{\parallel}$};
                    \draw[->] (2,-0.1) -- (2.5,-0.1) node[right] {\small $\nu_{\bot}$};
                \end{scope}
                \node at (0.5,0.5) {$\partial^* E$};
                \node at (3,1.7) {$\partial B(|x|)$};
                \node at (2.85,0.5) {\small $\nu^E(x)$};
        \end{tikzpicture}
        \end{center}
    \caption{Decomposition of $\nu^E$ as sum of its radial and tangential components.}
    \label{fig: decomp of nu^E}
\end{figure}
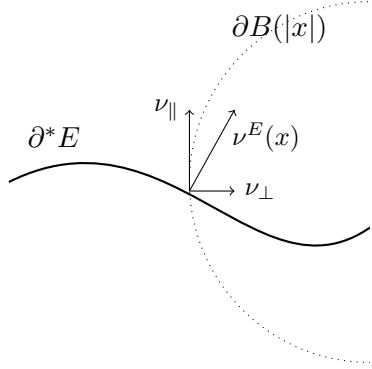

Let $E\subseteq \mathbb{R}^2$ be a measurable set; then we denote by $\Chi{E}$ its characteristic function and we say that $E$ is of finite perimeter if $\Chi{E}$ is in $BV(\mathbb{R}^2)$.
In this case, we call $P(E)=|D\Chi{E}|(\R^2)<+\infty$ its perimeter.
More in general, given a Borel set $A\subseteq \R^2$, the perimeter of $E$ in $A$ is defined as 
\[
P(E,A):=|D\Chi{E}|(A).
\]
We can give a different characterization of the relative perimeter of $E$ in $A$, by defining the {\it density points}.
Given $t\in[0,1]$, we denote by $E^t$ the set of points of density $t$ in $E$, namely
\[
E^t:=\left\{x\in\R^2\,\bigg\vert \lim_{\rho\to 0^+}\frac{\mathcal{H}^2(E\cap B(x,\rho))}{\pi \rho^2}=t\right\}.
\]
We define then the {\it essential boundary} of $E$ as $\partial ^eE:=E\setminus(E^1\cup E^0)$.
It turns out that, for any $A\subseteq\R^2$ Borel set,
\[
P(E,A):=\mathcal{H}^1(A\cap \partial^eE).
\]
For sets of finite perimeter, it is possible to define also the {\it reduced boundary} $\partial^*E$, that is the set of points in $\partial E$ such that there exists the {\it generalized outer unit normal} vector $\nu^E(x)$ given by
\[
\nu^E(x)=\lim_{\rho\to 0^+}\frac{D\Chi{E}(B(x,\rho))}{|D\Chi{E}|(B(x,\rho))}
\]
with $|\nu^E(x)|=1$;
in particular, it is well known that $\partial^*E\subseteq E^{\frac{1}{2}}$.
A fundamental result in this setting is given by the rectifiability theorem by De Giorgi (see \cite[Theorem 3.7]{A}), stating that, for all planar sets of finite perimeter, the reduced boudary is $1$-rectifiable and 
\[
D\Chi{E}=\nu^E\,\mathcal{H}^1\llcorner\partial^*E.
\]
The interested reader can find a complete presentation of finite-perimeter sets in \cite{A, AFP, M} and in the references therein.
Lastly, it will be useful to decompose the outer normal vector, for every $x\in\partial^*E$, as 
$$\nu^E(x)=\nu^E_\bot(x)+\nu^E_\parallel(x),$$ 
where $\nu^E_\bot(x)=(\nu^E(x)\cdot\hat{x})\hat{x}$ is the radial component of $\nu^E(x)$ along $\partial B(|x|)$ and $\nu^E_\parallel(x)$ is the tangential one ({\it cf.} Figure \ref{fig: decomp of nu^E}).
We remark that, if $\mathcal{H}^1(\partial B(r)\cap \{\nu^E_\parallel=0\})>0$, then the reduced boundary of $E$ coincides with $\partial B(r)$ in an arc of $\partial B(r)$ of positive measure.
In particular, this corresponds to a jump of the function $v$ at $r$.
\subsection*{Geometric measure theory results}
Some geometric measure theoretic results will be needed.
\begin{lemma}\label{lem partizione sup}
    Let $B\subseteq \mathbb{R}^n$ be a Borel set and let $\varphi_h,\varphi:B\to \mathbb{R}$, with $h\in \N$ be summable Borel functions such that $|\varphi_h|<|\varphi|$ for every $h$.
    Then it holds 
    \[
    \int_B \sup_h \varphi_h(x)dx=\sup_H \left\{\sum_{h\in H}\int_{A_h}\varphi_h(x)dx\right\},
    \]
    with the supremum made among all finite sets $H\subset\N$ and all finite partitions $\{A_h\}$, with $h\in H$, of $B$ in Borel sets.
\end{lemma} 

\begin{defin}
    For every $\varphi\in C_c(\R^n_0, \R^n)$, we decompose $\varphi=\varphi_{\bot}+\varphi_{\parallel}$ as the sum of its radial and tangential components, given by
    \[
    \varphi_{\bot}(x):=(\varphi (x)\cdot \hat{x})\hat{x},\hspace{1.5cm} \varphi_{\parallel}(x):=\varphi(x)-\varphi_{\bot}(x),
    \]
    where $\hat{x}=x/|x|$ denotes the direction of $x$. 
    Furthermore, we denote by div$_{\parallel}\varphi(x)$ the tangential divergence of $\varphi$ at $x$ along the sphere $\partial B(|x|)$, that is 
    \[
    {\rm div}_\parallel \varphi(x):={\rm div}\varphi(x)- \left(\nabla\varphi(x)\hat{x}\right)\cdot\hat{x},
    \] 
    where, with a little abuse of notation, we write $\nabla\varphi(x)\hat{x}$ to refer to $\nabla(\varphi(x)\cdot \hat{x})=D\varphi(x) \hat{x}$.
\end{defin}
By \cite[Lemma 3.2]{CPS}, in dimension $n=2$ the following expressions hold for every $x\in\R^2_0$:
\begin{equation}\label{formulae divergence}
    \bal
    {\rm div}\varphi_{\bot}(x)&=\left(\nabla\varphi(x)\hat{x}\right)\cdot\hat{x}+\frac{1}{|x|}\left(\varphi(x)\cdot\hat{x}\right)\\
    {\rm div}\varphi_{\parallel}(x)&={\rm div}_{\parallel}\varphi_{\parallel}(x).
    \eal
\end{equation}
We now consider the radial and tangential components of Radon measures as well, precisely in the case of $\mu=D\Chi{E}$, with $E\in \R^2$ a set of finite perimeter.
We set
\[
D_\bot\Chi{E}:=\nu^E_{\bot}\mathcal{H}^1\llcorner \partial^*E,\hspace{1.5cm} D_\parallel\Chi{E}:=\nu^E_{\parallel}\mathcal{H}^1\llcorner \partial^*E.
\]
By \cite[Lemma 3.4]{CPS}, for every $\varphi\in C^1_c(\R^2_0,\R^2)$ the following hold:
\begin{equation}\label{explicit DChi radial and tangential}
    \begin{aligned}
        \int_{\R^2_0}\varphi(x)dD_\parallel\Chi{E}&=-\int_{E}{\rm div}_\parallel\varphi_\parallel(x)dx\\
        \int_{\R^2_0}\varphi(x)dD_\bot\Chi{E}&=-\int_E\left(\nabla\varphi(x)\hat{x}\right)\cdot \hat{x} dx-\int_E\frac{\varphi(x)\cdot \hat{x}}{|x|} dx.
    \end{aligned}
\end{equation}

The next result is a special version of the coarea formula (see \cite[Theorem 2.93]{AFP}).
\begin{prop}\label{prop coarea parallel}
    Let $E$ be a set of finite perimeter in $\R^n$ and let $g:\R^n\to [0,+\infty]$ be a Borel function. 
    Then, writing $(\partial^*E)_r=\partial^*E\cap \partial B(r)$, it holds
    \[
    \int_{\partial^* E}g(x)|\nu_\parallel^E(x)|d\mathcal{H}^{n-1}(x)=\int_{0}^{+\infty}\left(\int_{(\partial^*E)_r}g(x)d\mathcal{H}^{n-2}(x)\right)dr.
    \]
\end{prop}

We close this section of preliminary results with a version of a result by Vol'pert (see \cite[Theorem 3.7]{CPS}):
\begin{theorem}\label{thm Vol'pert}
    Let $E\in \R^2$ be a set of finite perimeter and finite volume and let $v=v_E:(0,+\infty)\to [0,+\infty)$ be its circular distribution.
    Then, there exists a Borel subset of $\{\theta_v>0\}$, which is denoted by $G_E$, such that $\mathcal{H}^1(\{\theta_v>0\}\setminus G_E)=0$ and 
    \begin{itemize}
        \item For every $r\in G_E$, the slice $E_r$ is a set of finite perimeter in $\partial B(r)$;
        \item For every $r\in G_E$, the $0$-dimensional boundary of the slice $E_r$ coincides with the spherical slice of $\partial^*E$, namely, $\partial^* (E_r)=(\partial^*E)_r$ (hence we can write $\partial^*E_r$ without risk of ambiguity);
        \item For every $r\in G_E\cap \{0<\theta_v<\pi/2\}$ and for every $\omega\in \mathbb{S}^1$ such that $r\omega\in \partial^*E_r$, $|\nu_\parallel^E(r\omega)|>0$; moreover, $\nu_\parallel^E(r\omega)=\nu^{E_r}(r\omega)|\nu_\parallel^E(r\omega)|$, where $\nu^{E_r}(r\omega)$ is the direction corresponding to the angle $\vartheta(\omega)\pm\pi/2$.
    \end{itemize}
\end{theorem}

\section{Technical results}\label{sec. techn.}
We collect in this section some technical results that will be the core of the proof of Theorem \ref{thm sph rearr}.
The first one is a particular case of \cite[Lemma 4.1]{CPS}; for the sake of completeness, however, we show here the proof for the case $n=2$.
\begin{lemma}\label{lemma v in BV}
    Let $E\in \R^2$ be a set of finite perimeter and finite volume and let $v=v_E:(0,+\infty)\to [0,+\infty)$ be its circular distribution.
    Then $v$ is in $BV(0,+\infty)$ and $\theta_v$ is in $BV_{\rm loc}(0,+\infty)$.
    Moreover, for any $\psi:(0,+\infty)\to \R$ bounded Borel function, it holds 
    \begin{equation}\label{form rDthetav per psi}
        4\int_{0}^{+\infty}\psi(r)\cdot r dD\theta_v(r)=\int_{\R^2_0}\psi(|x|)\hat{x}dD_\bot\Chi{E}(x).
    \end{equation}
    In particular, for every $B\subseteq (0,+\infty)$ Borel set 
    \begin{equation}\label{dis |rDthetav| per B}
        |rD\theta_v|(B)\leq \frac{1}{4}|D\Chi{E}|(\phi(B\times\mathbb S^1)).
    \end{equation}
\end{lemma} 
\begin{proof}
    The proof is divided in two steps.
    \par {\bf Step 1.} We prove that $v$ belongs to $BV(0,+\infty)$.
    Trivially we have that $v$ is in $L^1(0,+\infty)$, due to the finiteness of the area of $E$. 
    Indeed, by definition of $v$ one has
    \[
    \int_0^{+\infty}v(r)dr=\int_0^{+\infty}\int_{\partial B(r)}\Chi{E}(x)d\mathcal{H}^1(x) dr=\mathcal{H}^2(E)<+\infty.
    \]
    Hence we have to prove that $Dv$ is a measure.
    Consider a function $\psi\in C^1_c(0,+\infty)$ with $|\psi|\leq 1$ and define 
    \[
    \varphi(x):=\psi(|x|)\hat{x}=\varphi_\bot(x).
    \]
    Then, thanks to \eqref{formulae divergence}, we have 
    \[
    \begin{aligned}
        {\rm div}\varphi(x)=&{\rm div}\varphi_\bot(x)
        =\big(\nabla\varphi(x)\hat{x}\big)\cdot\hat{x}+\frac{\varphi(x)\cdot\hat{x}}{|x|}
        =\big(\nabla(\psi(|x|)\hat{x})\hat{x}\big)\cdot\hat{x}+\frac{\psi(|x|)\hat{x}\cdot\hat{x}}{|x|}\\
        =&\left[\left(\psi'(|x|)\hat{x}\otimes \hat{x} +\psi(|x|)\frac{I-\hat{x}\otimes\hat{x}}{|x|}\right)\hat{x}\right]\cdot\hat{x}+\frac{\psi(|x|)}{|x|}\\
        =&\psi'(|x|)|\hat{x}|^4+\frac{\psi(|x|)}{|x|}(|\hat{x}|^2-|\hat{x}|^4)+\frac{\psi(|x|)}{|x|}=
        \psi'(|x|)+\frac{\psi(|x|)}{|x|}.
    \end{aligned}
    \]
    Hence, integrating against $\Chi{E}$ we obtain
    \[
    \int_{\R^2}\left(\psi'(|x|)+\frac{\psi(|x|)}{|x|}\right)\Chi{E}(x)dx=\int_{\R^2}{\rm div}(\psi(|x|)\hat{x})\Chi{E}(x)dx=-\int_{\R^2}\psi(|x|)\hat{x} dD_\bot\Chi{E}(x).
    \]
    In particular we have 
    \[
    \int_{\R^2}\psi'(|x|)\Chi{E}(x)dx=-\int_{\R^2}\psi(|x|)\hat{x}dD_\bot\Chi{E}(x)-\int_{\R^2}\frac{\psi(|x|)}{|x|}\Chi{E}(x)dx,
    \]
    where the left-hand side can be rewritten as 
    \[
    \begin{aligned}
        \int_{\R^2}\psi'(|x|)\Chi{E}(x)dx=&\int_{0}^{r}\psi'(r)\int_{\partial B(r)}\Chi{E}(x)d\mathcal{H}^1(x)\, dr=\int_{0}^{+\infty}\psi'(r)v(r)dr\\
        =&-\int_{0}^{+\infty}\psi(r)dDv(r).
    \end{aligned}
    \]
    Summarizing, we have the following identity:
    \[
    \int_{0}^{+\infty}\psi(r)dDv(r)=\int_{\R^2}\psi(|x|)\hat{x}dD_\bot\Chi{E}(x)+\int_{\R^2}\frac{\psi(|x|)}{|x|}\Chi{E}(x)dx.
    \]
    Now, the terms on the right-hand side can be estimated as 
    \[
    \begin{aligned}
        \int_{\R^2}\psi(|x|)\hat{x}dD_\bot\Chi{E}(x)=&\int_{\R^2}\varphi(x)dD_\bot\Chi{E}(x)\leq|D_\bot\Chi{E}|(\R^2)\leq P(E),\\
        \int_{\R^2}\frac{\psi(|x|)}{|x|}\Chi{E}(x)dx=&\int_{B(1)}\frac{\psi(|x|)}{|x|}\Chi{E}(x)dx+\int_{\R^2\setminus B(1)}\frac{\psi(|x|)}{|x|}\Chi{E}(x)dx\\
        \leq& \int_{0}^{1} \int_{0}^{2\pi}\frac{\psi(\rho)}{\rho}\Chi{E}(\rho\theta)\cdot \rho d\theta d\rho + \mathcal{H}^2(E)\leq 2\pi+\mathcal{H}^2(E).
    \end{aligned}
    \]
    So we found out that, for any $\psi\in C^1_c(0,+\infty)$ with $|\psi|\leq 1$, it holds
    \[
    \int_{0}^{+\infty}\psi(r)dDv(r)\leq 2\pi+\mathcal{H}^2(E)+P(E)<+\infty,
    \]
    thus obtaining that $v\in BV(0,+\infty)$.
    \par {\bf Step 2.} We conclude.
    First we notice that, since $v$ is in $BV(0,+\infty)$ and $r\mapsto 1/r$ is a smooth and locally bounded function in $(0,+\infty)$, we have that $\theta_v(r)=v(r)/4r$ is in $BV_{\rm loc}(0,+\infty)$.
    Moreover, thanks to Leibniz rule for $BV$ functions (see Proposition \ref{Leibniz rule BV}), we can write $Dv=4 D\big(\theta_v(r)\cdot r\big)$ as
    \[
    Dv=4\theta_v dr+4rD\theta_v.
    \]
    And so for any $\psi\in C^1_c(0,+\infty)$ with $|\psi|\leq 1$ we have
    \[
    \begin{aligned}
        \int_{0}^{+\infty}\psi(r)dDv(r)=&\int_{0}^{+\infty}\psi(r)\cdot \frac{v(r)}{r}dr+\int_{0}^{+\infty}4\psi(r)rdD\theta_v(r)\\
        =&\int_{\R^2}\psi(|x|)\frac{\psi(|x|)}{|x|}\Chi{E}(x)dx+\int_{0}^{+\infty}4\psi(r)rdD\theta_v(r).
    \end{aligned}
    \]
    This, combined with the computations in the first step, leads us to 
    \[
    \int_{0}^{+\infty}4\psi(r)rdD\theta_v(r)=\int_{0}^{+\infty}\psi(r)dDv(r)-\int_{\R^2}\psi(|x|)\frac{\psi(|x|)}{|x|}\Chi{E}(x)dx=\int_{\R^2}\psi(|x|)\hat{x}dD_\bot\Chi{E}(x),
    \]
    for all $\psi$ as above, and, by approximation, for all bounded Borel functions $\psi$, namely, \eqref{form rDthetav per psi} holds.
    Lastly, considering an open bounded set $B\subset(0,+\infty)$, for any test function $\psi\in C_c(B)$ with $|\psi|\leq 1$ it holds
    \[
    \int_B 4\psi(r)rdD\theta_v(r)=\int_{\phi(B\times \S^1)}\psi(|x|)\hat{x}dD_\bot\Chi{E}(x)\leq |D_\bot\Chi{E}|\big(\phi(B\times \S^1)\big).
    \]
    Then, taking the supremum on the left-hand side among all admissible $\psi$ and dividing both sides by $4$, we have 
    \[
    |rD\theta_v|(B)\leq \frac{1}{4} |D_\bot\Chi{E}|\big(\phi(B\times \S^1)\big),
    \]
    namely, \eqref{dis |rDthetav| per B} holds for open sets, and by approximation it holds for all  Borel sets in $(0,+\infty)$.
\end{proof}
In the next result we give more information about the measure $rD\eta_v$.
\begin{lemma}\label{lemma misura rDthetav}
    Let $E$ be as above. Then, for every Borel set $B$ in $(0,+\infty)$ we have
    \begin{equation}\label{eq: lemma misura primo claim}
        (4rD\theta_v)(B)=\int_{\partial^*E\cap \phi(B\times\mathbb S^1\cap \{\nu_\parallel^E=0\} )}\hat{x}\cdot \nu^E(x)d\mathcal{H}^1(x)+\int_B dr \int_{\partial^*E_r\cap\{\nu_\parallel^E\neq 0\}}\frac{\hat{x}\cdot \nu^E(x)}{|\nu_\parallel^E(x)|}d\mathcal{H}^0(x).
    \end{equation}
    In addition, the restriction of $rD\theta_v$ to the set $G_{F_v}$ gives
    \begin{equation}\label{eq: lemma misura secondo claim}
        rD\theta_v\llcorner G_{F_v}=r\theta_v'(r)\mathcal{L}^1\llcorner G_{F_v}=\frac{1}{4}\,\mathcal{H}^0(S_{\theta_v(r)}(r)) \frac{\hat{x}\cdot \nu^{F_v}(x)}{|\nu_\parallel^{F_v}(x)|}\mathcal{L}^1\llcorner G_{F_v},
    \end{equation}
    for any $x\in (\partial^*F_v)_r$.
\end{lemma}
\begin{proof}
    Take a Borel set $B\subseteq (0,+\infty)$ and let $\psi$ be the Borel map defined by $\psi(r):=\Chi{B}(r)$.
    Then, by \eqref{form rDthetav per psi} and by definition of $D_\bot\Chi{E}$ we have
    \[
    \begin{aligned}
        (4rD\theta_v)(B)&=\int_{0}^{+\infty}4\Chi{B}(r)rdD\theta_v(r)
        =\int_{\R ^2_0}\Chi{B}(|x|)\hat{x}dD_\bot\Chi{E}(x)
        =\int_{\partial^*E}\Chi{B}(|x|)\hat{x}\cdot\nu^E_\bot d\mathcal H^1(x)\\
        &=\int_{\partial^*E\cap\phi(B\times \mathbb S^1)}\hat{x}\cdot\nu^E_\bot d\mathcal H^1(x)
        =\int_{\partial^*E\cap\phi(B\times \mathbb S^1)}\hat{x}\cdot\nu^E d\mathcal H^1(x)\\
        &=\int_{\partial^*E\cap\phi(B\times \mathbb S^1)\cap \{\nu_\parallel^E=0\}}\hat{x}\cdot\nu^E d\mathcal H^1(x)+\int_{\partial^*E\cap\phi(B\times \mathbb S^1)\cap \{\nu_\parallel^E\neq0\}}\hat{x}\cdot\nu^E d\mathcal H^1(x).
    \end{aligned}
    \]
    The second term of the last sum can be exploited via the coarea formula in Proposition \ref{prop coarea parallel} as follows:
    \[
    \begin{aligned}
        \int_{\partial^*E\cap\phi(B\times \mathbb S^1)\cap \{\nu_\parallel^E\neq0\}}\hat{x}\cdot\nu^E d\mathcal H^1(x)
        &=\int_{0}^{+\infty}dr \int_{\partial^*E_r\cap\phi(B\times \mathbb S^1)\cap \{\nu_\parallel^E\neq0\}}\frac{\hat{x}\cdot\nu^E}{|\nu_\parallel^E(x)|} d\mathcal H^0(x)\\
        &=\int_B dr \int_{\partial^*E_r\cap \{\nu_\parallel^E\neq0\}}\frac{\hat{x}\cdot\nu^E}{|\nu_\parallel^E(x)|} d\mathcal H^0(x),
    \end{aligned}
    \] 
    thus obtaining \eqref{eq: lemma misura primo claim}.
    Concerning \eqref{eq: lemma misura secondo claim}, we notice that, for any $r\in G_{F_v}$ such that $\theta_v(r)\in \left\{0,\pi/2\right\}$ it holds $\mathcal H^0((\partial^*F_v)_r)=0$,
    while for all $r\in G_{F_v}\cap \{0<\theta_v(r)<\pi/2\}$ we have that 
    \[
    \frac{\hat{x}\cdot \nu^E(x)}{|\nu_\parallel^E(x)|}
    \]
    is constant on $(\partial^*F_v)_r$ by symmetry of $F_v$.
    Hence we have that 
    \[
    (4rD\theta_v)\llcorner G_{F_v}=\left(\int_{(\partial^*F_v)_r}\frac{\hat{x}\cdot\nu^{F_v}}{|\nu_\parallel^{F_v}(x)|} d\mathcal H^0(x)\right)\,\mathcal{L}^1\llcorner G_{F_v}=\mathcal H^0((\partial^*F_v)_r)\frac{\hat{x}\cdot\nu^{F_v}}{|\nu_\parallel^{F_v}(x)|}\,\mathcal{L}^1\llcorner G_{F_v}
    \]
    for any $x\in (\partial^*F_v)_r$.
\end{proof}
The last technical result that we include gives us a first estimate for the local perimeter of the rearrangement $F_v$.
\begin{prop}\label{prop estim P(F)}
    Let $E,v$ be as above. Then, $F_v$ is a set of finite perimeter in $\R^2$ and for every Borel set $B$ in $(0,+\infty)$ the following inequality holds:
    \begin{equation}\label{ineq estim P(F)}
        P(F_v, \phi(B\times\mathbb{S}^1))\leq |4rD\theta_v|(B)+|D_\parallel\Chi{F_v}|(\phi(B\times \mathbb S^1)).
    \end{equation}   
\end{prop}
\begin{proof}
    By Lemma \ref{lemma v in BV} we know that $v$ is in $BV(0,+\infty)$, hence we can find a sequence of non-negative functions $\{v_j\}_j\subseteq C^{\infty}_c(0,+\infty)$ such that
    \[v_j\rightarrow v\quad \mathcal{H}^1-a.e.,\hspace{1.5cm} |Dv_j|\overset{*}{\rightharpoonup}|Dv|.\]
    Consider an open set $\Omega\subseteq (0,+\infty)$ and a test function $\varphi\in C^1_c(\phi(\Omega\times\mathbb S^1);\R^2)$ with $\|\varphi\|_\infty\leq 1$.
    For every $j$ we have, using \eqref{explicit DChi radial and tangential}
    \[
    \begin{aligned}
        &\int_{\phi(\Omega\times \S^1)}\Chi{F_{v_j}}(x){\rm div}\varphi(x)dx=
        -\int_{\phi(\Omega\times \S^1)}\varphi(x)dD\Chi{F_{v_j}}(x)\\
        =&\int_{\phi(\Omega\times \S^1)}\Chi{F_{v_j}}{\rm div}_\parallel\varphi_\parallel(x)dx+\int_{\phi(\Omega\times \S^1)}\Chi{F_{v_j}}(\nabla\varphi(x)\hat{x})\cdot \hat{x}dx+\int_{\phi(\Omega\times \S^1)}\Chi{F_{v_j}}\frac{\varphi(x)\cdot\hat{x}}{|x|}dx.
    \end{aligned}
    \]
    Our aim is to estimate the three terms above and then take the limit as $j\to +\infty$.
    Let us define the function $V_j(r)$ as 
    \[
    V_j(r)=\int_{D_{\theta_{v_j}(r)}(r)}\varphi(x)\cdot\hat{x}d\mathcal{H}^1(x)=\int_{D_{\theta_{v_j}(r)}(1)}r\varphi(r\omega)\cdot\hat{\omega}d\mathcal{H}^1(\omega).
    \]
    The rest of the proof is divided into several steps.
    \par {\bf Step 1.} We prove that $V_j$ is Lipschitz and has compact support.
    We start by noticing that the support of $V_j$ is contained in
    \[
    \Lambda supp(\varphi):=\{r\in (0,+\infty)\,:\,supp(\varphi)\cap \partial B(r)\neq \emptyset\},
    \]
    which is compact by compactness of the support of $\varphi$.
    Now consider $r_1,r_2$ in $supp(V_j)$ and assume (without loss of generality) that $\theta_{v_j}(r_1)\geq \theta_{v_j}(r_2)$; then it holds
    \[
    \begin{aligned}
        |V_j(r_1)-V_j(r_2)|\leq
        &\,\int_{D_{\theta_{v_j}(r_1)}(1)} \big|r_1 \varphi(r_1\omega)\cdot \omega - r_2 \varphi(r_2\omega)\cdot \omega\big|d\mathcal{H}^1(\omega) \\
        &\,+r_2\left|\int_{D_{\theta_{v_j}(r_2)}(1)}\varphi(r_2\omega)\cdot \omega d\mathcal{H}^1(\omega) - \int_{D_{\theta_{v_j}(r_1)}(1)}\varphi(r_2\omega)\cdot\omega d\mathcal H^1(\omega)\right|\\
        \leq\, & c|r_1-r_2| + r_2\int_{D_{\theta_{v_j}(r_1)}(1)\setminus D_{\theta_{v_j}(r_2)}(1)}|\varphi(r_2\omega)\cdot \omega|d\mathcal{H}^1(\omega)\\
        \leq\, & c|r_1-r_2| + r_2 \mathcal{H}^1 \big(D_{\theta_{v_j}(r_1)}(1)\setminus D_{\theta_{v_j}(r_2)}(1)\big)\\
        = &\, c|r_1-r_2|+4r_2|\theta_{v_j}(r_1)-\theta_{v_j}(r_2) |\\
        \leq&\, c'|r_1-r_2|,
    \end{aligned}
    \]
    where we used that, due to the compactness of $supp(\theta_{v_j})$, the Lipschitzianity of $v_j$ implies Lipschitzianity of $\theta_{v_j}$ for all $j$'s.
    This concludes the first step.
    \par {\bf Step 2.} We find now an expression for $V_j'$.
    In the previous step we proved that $V_j$ is Lipschitz and, in particular, this gives differentiability $\mathcal{H}^1$-$a.e.$. 
    Moreover, since the $v_j$'s are $C^{\infty}$ functions, then $\theta_{v_j}=v_j/4r \in C^{\infty}((0,+\infty))$ and we can compute, keeping in mind that $S_{\sigma}(1)$ denotes the four points $(\pm \cos \sigma, \pm \sin \sigma)$,
    \[
    \begin{aligned}
        V_j'(r)=&\frac{d}{dr}\left(r\int_{0}^{\theta_{v_j}(r)}d\sigma\int_{S_\sigma(1)}\varphi(r\omega)\cdot \omega d\mathcal{H}^0(\omega)\right)
        =\int_{0}^{\theta_{v_j}(r)}d\sigma\int_{S_\sigma(1)}\varphi(r\omega)\cdot \omega d\mathcal{H}^0(\omega)\\
        &+r\theta_{v_j}'(r)\int_{S_{\theta_{v_j}(r)}(1)}\varphi(r\omega)\cdot \omega d\mathcal{H}^0(\omega)
        +r\int_{0}^{\theta_{v_j}(r)}d\sigma\int_{S_\sigma(1)}(\nabla\varphi(r\omega)\omega)\cdot \omega d\mathcal{H}^0(\omega)\\
        =&\int_{D_{\theta_{v_j}(r)}(1)}\varphi(r\omega)\cdot \omega d\mathcal{H}^1(\omega)
        +r\theta_{v_j}'(r)\int_{S_{\theta_{v_j}(r)}(1)}\varphi(r\omega)\cdot \omega d\mathcal{H}^0(\omega)\\
        &+r\int_{D_{\theta_{v_j}(r)}(1)}(\nabla\varphi(r\omega)\omega)\cdot \omega d\mathcal{H}^1(\omega).
    \end{aligned}
    \]
    Summarizing, we showed that for $\mathcal{H}^1$-$a.e.$ $r$ 
    \begin{equation}\label{expr Vj'}
        \begin{aligned}
            V_j'(r)=&\int_{D_{\theta_{v_j}(r)}(1)}\varphi(r\omega)\cdot \omega d\mathcal{H}^1(\omega)
            +r\theta_{v_j}'(r)\int_{S_{\theta_{v_j}(r)}(1)}\varphi(r\omega)\cdot \omega d\mathcal{H}^0(\omega)\\
            &+r\int_{D_{\theta_{v_j}(r)}(1)}(\nabla\varphi(r\omega)\omega)\cdot \omega d\mathcal{H}^1(\omega).
        \end{aligned}
    \end{equation}
    \par {\bf Step 3.} We prove that
    \begin{equation}\label{eq a+b=-c}
        \begin{aligned}
            &\int_{\phi(\Omega\times \S^1)}\Chi{F_{v_j}}\big(\nabla\varphi(x)\hat{x}\big)\cdot \hat{x}dx+\int_{\phi(\Omega\times \S^1)}\Chi{F_{v_j}}\frac{\varphi(x)\cdot\hat{x}}{|x|}dx\\
            =&-\int_{\Omega}r\theta_{v_j}'(r)\left(\int_{S_{\theta_{v_j}(r)}(1)}\varphi (r\omega)\cdot \omega d\mathcal{H}^0(\omega)\right) dr.
        \end{aligned}
    \end{equation}
    We integrate \eqref{expr Vj'} over $\Omega$ and, using the fact that $V_j$ has compact support, we obtain
    \begin{equation}\label{eq Vj' integr}
        \begin{aligned}
            0=&\int_{\Omega}dr\int_{D_{\theta_{v_j}(r)}(1)}\varphi(r\omega)\cdot \omega d\mathcal{H}^1(\omega)
            +\int_{\Omega}r\theta_{v_j}'(r)\left(\int_{S_{\theta_{v_j}(r)}(1)}\varphi(r\omega)\cdot \omega d\mathcal{H}^0(\omega)\right) dr\\
            &+\int_{\Omega}r\left(\int_{D_{\theta_{v_j}(r)}(1)}(\nabla\varphi(r\omega)\omega)\cdot \omega d\mathcal{H}^1(\omega)\right) dr.
        \end{aligned}
    \end{equation}
    Now, the first and the third term can be rewritten as 
    \[
    \begin{aligned}
        \int_{\Omega}dr\int_{D_{\theta_{v_j}(r)}(1)}\varphi(r\omega)\cdot \omega d\mathcal{H}^1(\omega)&=\int_{\phi(\Omega\times \S^1)}\Chi{F_{v_j}}(x)\frac{\varphi(x)\cdot \hat{x}}{|x|}dx,\\
        \int_{\Omega}r\left(\int_{D_{\theta_{v_j}(r)}(1)}(\nabla\varphi(r\omega)\omega)\cdot \omega d\mathcal{H}^1(\omega)\right) dr&=\int_{\phi(\Omega\times \S^1)}\Chi{F_{v_j}}(x)(\nabla\varphi(x)\hat{x})\cdot\hat{x}dx,
    \end{aligned}
    \]
    thus, using these two expressions in \eqref{eq Vj' integr}, we have the claim.
    \par {\bf Step 4.}
    We prove that 
    \begin{equation}\label{ineq step4}
        \int_{\phi (\Omega\times \S ^1)}\Chi{F_{v_j}}(x){\rm div}\varphi(x)dx \leq \int_{\Omega}\mathcal{H}^0\big(S_{\theta_{v_j}(r)}(r)\big)dr + |4rD\theta_{v_j}|(\Lambda supp\varphi),
    \end{equation}
    where $\Lambda supp\varphi$ is the set defined in the first step.\\
    We know from the previous results and from \eqref{eq a+b=-c} that 
    \begin{equation*}
        \begin{aligned}
            &\int_{\phi (\Omega\times \S ^1)}\Chi{F_{v_j}}(x){\rm div}\varphi(x)dx\\
            =&\int_{\phi(\Omega\times \S^1)}\Chi{F_{v_j}}{\rm div}_\parallel\varphi_\parallel(x)dx+\int_{\phi(\Omega\times \S^1)}\Chi{F_{v_j}}(\nabla\varphi(x)\hat{x})\cdot \hat{x}dx+\int_{\phi(\Omega\times \S^1)}\Chi{F_{v_j}}\frac{\varphi(x)\cdot\hat{x}}{|x|}dx\\
            =&\int_{\phi(\Omega\times \S^1)}\Chi{F_{v_j}}{\rm div}_\parallel\varphi_\parallel(x)dx - \int_{\Omega}r\int_{S_{\theta_{v_j}(r)}(1)}\theta_{v_j}'(r)\varphi(r\omega)\cdot\omega d\mathcal{H}^0(\omega)\, dr.
        \end{aligned}
    \end{equation*}
Now we notice that, whenever $supp(\varphi)\cap \partial B(r)\neq \emptyset$, we have 
\[
\int_{S_{\theta_{v_j}(r)}(1)}\varphi(r\omega)\cdot \omega d\mathcal{H}^0(\omega)\leq \mathcal{H}^0(S_{\theta_{v_j}(r)}(1)) \leq 4;
\]
hence, we can write
\[
-\int_{\Omega} r\theta_{v_j}'(r)\int_{S_{\theta_{v_j}(r)}(1)}\varphi(r\omega)\cdot \omega d\mathcal{H}^0(\omega)\,dr\leq\int_{\Lambda supp\varphi}4r\cdot |\theta_{v_j}'(r)|dr \leq |4rD\theta_{v_j}|(\Lambda supp\varphi).
\]
On the other hand, regarding the integral of the tangential divergence, we can apply the divergence theorem on the circle and we have
\begin{equation*}
    \begin{aligned}
        \int_{\phi(\Omega\times \S^1)}\Chi{F_{v_j}}{\rm div}_\parallel\varphi_\parallel(x)dx &=\int_{\Omega}\int_{D_{\theta_{v_j}(r)}(r)} {\rm div}_\parallel \varphi_\parallel(x)d\mathcal{H}^1(x)\, dr \\
        &=\int_{\Omega} \int_{S_{\theta_{v_j}(r)}(r)}\varphi_\parallel(x)\cdot \nu^{*}(x) d\mathcal{H}^0(x)\,dr
        \leq \int_{\Omega}\mathcal{H}^0(S_{\theta_{v_j}(r)}(r))dr,
    \end{aligned}
\end{equation*} 
where $\nu^*$ denotes the normal outward vector to $D_{\theta_{v_j}(r)}(r)$ at $x$.
Putting these two estimates together we obtain \eqref{ineq step4}, as we wanted.
\par {\bf Step 5.} We prove that $F_{v}$ has locally-finite perimeter.
By definition of the sequence $\{v_j\}$, we have that $\Chi{F_{v_j}}\to \Chi{F_v}$ $\mathcal{H}^2-a.e.$, as well as $\theta_{v_j}\to\theta_v$ $\mathcal{H}^1-a.e.$ and $|rD\theta_{v_j}|\overset{*}{\rightharpoonup}|rD\theta_v|$.
Hence we have, by Step 4,
\[
\begin{aligned}
    \int_{\phi(\Omega\times\S ^1)}\Chi{F_v}(x){\rm div}\varphi(x)dx&
    =\limsup_j \int_{\phi(\Omega\times\S ^1)}\Chi{F_{v_j}}(x){\rm div}\varphi(x)dx\\
    &\leq \limsup_j \int_{\Omega}\mathcal{H}^0(S_{\theta_{v_j}(r)}(r))dr + \limsup_j |4rD\theta_{v_j}|(\Lambda supp\varphi).
\end{aligned}
\]
Let us focus on the first term, where we have
\[
\limsup_j \int_{\Omega}\mathcal{H}^0(S_{\theta_{v_j}(r)}(r))dr = \int_{\Omega}\mathcal{H}^0(S_{\theta_{v}(r)}(r))dr \leq 4 \mathcal{H}^1\left(\Omega \cap \left\{0<\theta_v< \frac{\pi}{2}\right\}\right)\leq C
\]
where, in the last inequality, we used the hypothesis that $E$ has finite volume and finite perimeter.
Regarding the second term we have, by compactness of $\Lambda supp \varphi$ and using \eqref{dis |rDthetav| per B},
\[
\limsup_j |rD\theta_{v_j}|(\Lambda supp\varphi)=|rD\theta_{v}|(\Lambda supp\varphi)\leq |rD\theta_{v}|(\Omega)\leq \frac{1}{4}|D_\bot\Chi{E}|(\phi(\Omega\times \S^1)).
\] 
Hence, by taking the supremum among all possible test functions $\varphi$ we found that 
\[
P\big(F_v; \phi(\Omega\times \S ^1)\big)\leq P\big(E;\phi(\Omega\times \S ^1)\big)+C,
\]
which tells us that $F_v$ has locally-finite perimeter.
\par {\bf Step 6.} We conclude, proving \eqref{ineq estim P(F)}.
Consider again a test function $\varphi$ as above.
For any $j$ we have 
\[
\int_{\phi(\Omega\times \S^1)}\Chi{F_{v_j}}(x){\rm div}\varphi(x) dx\leq \int_{\phi(\Omega\times \S^1)}\Chi{F_{v_j}}(x){\rm div}_\parallel\varphi_\parallel(x) dx + |4rD\theta_{v_j}|(\Lambda supp\varphi),
\] 
hence taking the superior limit as $j\to +\infty$ we have
\[
\int_{\phi(\Omega\times \S^1)}\Chi{F_{v}}(x){\rm div}\varphi(x) dx\leq \big|D_\parallel\Chi{F_{v}}\big|(\phi(\Omega\times \S^1)) + |4rD\theta_{v}|(\Omega).
\]
Now, taking the supremum on the left-hand side among all possible test functions we have \eqref{ineq estim P(F)} for open subsets of $(0,+\infty)$.
If $B$ is a Borel subset of $(0,+\infty)$, by regularity we can take the infimum of \eqref{ineq estim P(F)} over all open $\Omega$ containing $B$, and the inequality holds true for $B$ as well.
This concludes the proof of the proposition.
\end{proof}

\section{Proof of Theorem \ref{thm sph rearr}}\label{sec. proof thm}
We are finally ready to prove the main result of this paper.
Consider a Borel set $B\subseteq (0,+\infty)$.
We consider two cases, depending on the relation between $B$ and the ``good set'' $G_{F_v}$ in Vol'pert's Theorem \ref{thm Vol'pert}.
\par {\bf Case 1.} Let us first consider the case where $B\cap G_{F_v}=\emptyset$.
By Proposition \ref{prop coarea parallel} we know that 
\[
\begin{aligned}
    |D_\parallel \Chi{F_v}|(\phi(B\times \S ^1))
    &=\int_{\phi(B\times \S ^1)}|\nu_\parallel^{F_v}(x)|d(\mathcal{H}^1\llcorner\partial^*F_v)(x)
    =\int_{\phi(B\times \S ^1)\cap \partial^*F_v}|\nu_\parallel^{F_v}(x)|d\mathcal{H}^1(x)\\
    &=\int_{B}dr\int_{(\partial^*F_v)_r}1  d\mathcal{H}^0(x)
    =\int_B\mathcal{H}^0\big((\partial^*F_v)_r\big)dr.
\end{aligned}
\]
Now we just notice that for almost every $r\in B$ it must hold $\theta_v(r)=0$ (and consequently $(\partial^*F_v)_r=\emptyset$) since $\mathcal{H}^1(\{\theta_v>0\}\setminus G_{F_v})=0$, and so we have 
\[
|D_\parallel \Chi{F_v}|(\phi(B\times \S ^1))=0.
\]
Then, combining this with Proposition \ref{prop estim P(F)} and \eqref{dis |rDthetav| per B} we have
\[
P(F_v;\phi(B\times \S ^1))\leq |4rD\theta_v|(B)\leq |D_\bot \Chi{E}|(\phi(B\times \S ^1))\leq P(E;\phi(B\times \S ^1)).
\]
\par {\bf Case 2.} Now we consider the case where $B\subseteq G_{F_v}$; the general case can then be managed by decomposition.
This part of the proof is divided in steps.
\par {\bf Step 1.} First, we see that 
\begin{equation}\label{case 2, step 1}
    P(E,\phi(B\times \S ^1))\geq P\big(E;\phi(B\times \S ^1)\cap \{\nu_\parallel^E=0\}\big)+\int_B \sqrt{g_E^2(r)+p_E^2(r)}dr,
\end{equation}
where $p_E$ and $g_E$ are defined as follows
\begin{equation}\label{def pE, gE}
    p_E(r):=\mathcal{H}^0(\partial^*E_r),\hspace{1.5cm} 
    g_E(r):=\int_{\partial^*E_r}\frac{\hat{x}\cdot \nu^E(x)}{|\nu_\parallel^E(x)|}d\mathcal{H}^0(x).
\end{equation}
We start by decomposing 
\[
P\big(E;\phi(B\times\S^1)\big)=P\big(E;\phi(B\times\S^1)\cap\{\nu_\parallel^E=0\}\big)+P\big(E;\phi(B\times\S^1)\cap\{\nu_\parallel^E\neq 0\}\big);
\]
now we focus on the second term, which by the coarea formula becomes
\[
\begin{aligned}
    P(E;\phi(B\times\S^1)\cap\{\nu_\parallel^E\neq 0\})
    &=\int_Bdr\int_{\partial^*E_r}\frac{1}{|\nu_\parallel^E(x)|}d\mathcal{H}^0(x)\\
    &=\int_Bdr\int_{\partial^*E_r}\sqrt{1+\left(\frac{\hat{x}\cdot \nu^E(x)}{|\nu_\parallel^E(x)|}\right)^2}d\mathcal{H}^0(x),
\end{aligned}
\]
where in the last equality we used that $1=|\nu_\parallel^E|^2+|\nu_\bot^E|^2$. 
Now, being $t\mapsto \sqrt{1+t^2}$ a convex function, we can apply Jensen's inequality and get to 
\[
\begin{aligned}
    \int_Bdr\int_{\partial^*E_r}\sqrt{1+\left(\frac{\hat{x}\cdot \nu^E(x)}{|\nu_\parallel^E(x)|}\right)^2}d\mathcal{H}^0
    &\geq \int_B\mathcal{H}^0(\partial^*E_r)\sqrt{1+\left(\frac{1}{\mathcal{H}^0(\partial^*E_r)}\cdot \int_{\partial^*E_r}\frac{\hat{x}\cdot\nu^E(x)}{|\nu_\parallel^E(x)|}\right)^2}\, dr\\
    &=\int_B p_E(r)\sqrt{1+\frac{g^2_E(r)}{p_E^2(r)}}dr
    =\int_B\sqrt{p_E^2(r)+g_E^2(r)}dr,
\end{aligned}
\]
as we wanted.
\par {\bf Step 2.} We prove that 
\begin{equation}\label{stima P(E)}
    P(E;\phi(B\times\S^1))\geq \int_B\sqrt{p_E^2(r)+(4r\theta'_v(r))^2}\, dr.
\end{equation}
Let $\{A_h\}_{h\in H}$ be a finite Borel partition of $B$, with $A_h\subseteq G_{F_v}$ for every $h$.
Then by Lemma \ref{lemma misura rDthetav} we know that for every $h\in H$
\[ rD\theta_v\llcorner A_h=r\theta_v'\mathcal{L}^1 \llcorner A_h.\]
In particular, for any $\{a_h\}\subseteq\R$ it holds
\[
\begin{aligned}    
    &4\int_{A_h}a_hr\theta'_v(r)dr=4\int_{A_h}a_hrdD\theta_v(r)\\
    =&\int_{\partial^*E\cap\phi(A_h\times\S^1)\cap \{\nu^E_\parallel=0\}}a_h \hat{x}\cdot \nu^E(x)d\mathcal{H}^1
    +\int_{A_h}dr\int_{\partial^*E_r\cap\phi \{\nu^E_\parallel\neq0\}}a_h\, \frac{\hat{x}\cdot \nu^E(x)}{|\nu^E_\parallel(x)|}d\mathcal{H}^0;
\end{aligned}
\]
hence we have, reminding the definition of $g_E$ given in \eqref{def pE, gE},
\[
\int_{A_h}4a_hr\theta'_v(r)dr=\int_{\partial^*E\cap\phi(A_h\times\S^1)\cap \{\nu^E_\parallel=0\}}a_h \hat{x}\cdot \nu^E(x)d\mathcal{H}^1+\int_{A_h}a_h g_E(r)dr.
\]
Before going on, we notice that in $\{\nu^E_\parallel=0\}$, it holds $\hat{x}\cdot \nu^E(x)=1$, hence
\[
\int_{\partial^*E\cap\phi(A_h\times\S^1)\cap \{\nu^E_\parallel=0\}}a_h \hat{x}\cdot \nu^E(x)d\mathcal{H}^1=a_hP\big(E;\phi (A_h\times \S^1)\cap \{\nu^E_\parallel=0\}\big).
\]
Moreover, a simple calculation shows that
\begin{equation}\label{convex dual}
    \sqrt{1+t^2}=\sup_{h\in\N}\left\{a_h t+\sqrt{1-a_h^2}\right\},
\end{equation}
among a dense set $\{a_h\}$ in $(-1,1)$.
Thus we have
\[
\begin{aligned}
    &\sum_{h\in H}\int_{A_h} 4a_hr\theta'_v(r)dr+\int_{A_h}p_E(r)\sqrt{1-a_h^2}dr\\
    =&\sum_{h\in H}\int_{\partial^*E\cap\phi(A_h\times\S^1)\cap \{\nu^E_\parallel=0\}}a_h \hat{x}\cdot \nu^E(x)d\mathcal{H}^1+\int_{A_h}a_h g_E(r)+p_E(r)\sqrt{1-a_h^2}dr\\
    =&\sum_{h\in H} a_hP\big(E;\phi (A_h\times \S^1)\cap \{\nu^E_\parallel=0\}\big)+\int_{A_h}p_E(r)\left(a_h\frac{g_E(r)}{p_E(r)}+\sqrt{1-a_h^2}\right)dr\\
    \leq& \sum_{h\in H} a_hP\big(E;\phi (A_h\times \S^1)\cap \{\nu^E_\parallel=0\}\big)+\int_{A_h}\sqrt{p_E^2(r)+g_E^2(r)}dr\\
    \leq& P\big(E;\phi (B\times \S^1)\cap \{\nu^E_\parallel=0\}\big)+\int_{B}\sqrt{p_E^2(r)+g_E^2(r)}dr.
\end{aligned}
\]
Hence, taking the supremum among all possible finite partitions of $B$ we have 
\[
\begin{aligned}
    &\sup_{H\subseteq \N}\left(\sum_{h\in H}\int_{A_h}4 a_hr\theta'_v(r)dr+\int_{A_h}p_E(r)\sqrt{1-a_h^2}dr\right)\\
    \leq& P\big(E;\phi (B\times \S^1)\cap \{\nu^E_\parallel=0\}\big)+\int_{B}\sqrt{p_E^2(r)+g_E^2(r)}dr 
    \leq P\big(E;\phi(B\times \S^1)\big),
\end{aligned}
\]
where the last inequality is \eqref{case 2, step 1}.
On the other hand, the left-hand side can be rewritten, thanks to Lemma \ref{lem partizione sup}, as 
\[\begin{aligned}
    &\sup_{H\subseteq \N}\left(\sum_{h\in H}\int_{A_h} 4a_hr\theta'_v(r)dr+\int_{A_h}p_E(r)\sqrt{1-a_h^2}dr\right)
    \\=&\int_B \sup_h\left(4a_hr\theta_v'(r)+p_E(r)\sqrt{1-a_h^2}\right)dr=\int_B\sqrt{p_E^2(r)+(4r\theta'_v(r))^2}dr
\end{aligned}\]
where in the last equality we used again \eqref{convex dual}.
This gives us the claim of this step.
\par
\begin{figure}[!ht]
    \begin{center}
        \begin{tikzpicture}
            \draw (0.9848,0.1736) arc[start angle=10, end angle=175, radius=1];
            \draw (0.9848,-0.1736) arc[start angle=350, end angle=185, radius=1];
            \draw (0.9848,0.1736) -- (2.9848,0.1736);
            \draw (2.9848,-0.1736) -- (0.9848,-0.1736);
            \draw (-0.9962,0.0871) -- (-2.9962,0.0871) -- (-2.9962,-0.0871) -- (-0.9962,-0.0871);
            \draw (2.9848,-0.1736) arc[start angle=-5, end angle=5, radius=1.9918];
            \node at (3.3, 0.05) {$2r$};
            \node at (-3.4, 0.05) {$2\alpha r$};
            \node at (2, 0.5) {$l$};
            \node at (-2, 0.5) {$l$};
        \end{tikzpicture}
    \end{center}
    \vspace{0.3cm}
    \begin{center}
        \begin{tikzpicture}
            \draw (0.9914,0.1305) arc[start angle=7.5, end angle=172.5, radius=1];
            \draw (0.9914,-0.1305) arc[start angle=352.5, end angle=187.5, radius=1];
            \draw (0.9914,0.1305) -- (2.9914, 0.1305) -- 
            (2.9914,-0.1305) -- (0.9914,-0.1305);
            \draw (-0.9914,0.1305) -- (-2.9914, 0.1305) -- 
            (-2.9914,-0.1305) -- (-0.9914,-0.1305);
            \node at (3.3, 0.05) {$2\tilde{r}$};
            \node at (-3.3, 0.05) {$2\tilde{r}$};
            \node at (2, 0.5) {$l$};
            \node at (-2, 0.5) {$l$};
        \end{tikzpicture}
    \end{center}
    \caption{The example showing that, in general, the double spherical cap symmetrization does not decrease the perimeter in dimension higher than $2$.}
    \label{fig counterex double sph high dim}
\end{figure}
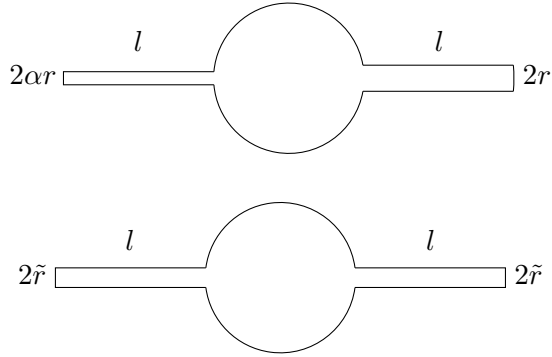

 {\bf Step 3.} We prove \eqref{strong ineq sph rearr}, {\it i.e.}
\[
P(F_v; \phi(B\times \S^1))\leq P(E; \phi(B\times \S^1)) + 2\mathcal{H}^1(\Gamma_E\cap B),
\]
where we remind that $\Gamma_E$ is defined as
\[
\Gamma_E=\{r\in (0,+\infty)\,|\, E_r \text{ is connected,  } 0<\mathcal{H}^1(E_r)<2\pi r\}.
\]
Let us show that \eqref{stima P(E)} is an equality if we choose $E=F_v$, that is,
\[
P(F_v; \phi(B\times \S^1))=\int_B\sqrt{p_{F_v}^2(r)+(4r\theta'_v(r))^2}dr.
\]
Indeed 
\[
\begin{aligned}
    P(F_v; \phi(B\times \S^1))&=\mathcal{H}^1(\partial^*F_v\cap \phi(B\times \S^1))=
    \int_{B\cap \{0<\theta_v<\pi/2\}}dr\int_{(\partial^*F_v)_r}\frac{1}{|\nu_\parallel^{F_v}(x)|}d\mathcal{H}^0(x)
    \\&=\int_{B\cap \{0<\theta_v<\pi/2\}}dr\int_{(\partial^*F_v)_r}\sqrt{1+\frac{|\nu_\bot^{F_v}(x)|^2}{|\nu_\parallel^{F_v}(x)|^2}}d\mathcal{H}^0(x).
\end{aligned}
\]
Now, by Lemma \ref{lemma misura rDthetav}, for all $x\in(\partial^*F_v)_r$ it holds
\[
\frac{|\nu_\bot^{F_v}(x)|^2}{|\nu_\parallel^{F_v}(x)|^2}=\left(\frac{4r\theta_v'(r)}{\mathcal{H}^0(S_{\theta_v(r)}(r))}\right)^2=\left(\frac{4r\theta_v'(r)}{p_{F_v}(r)}\right)^2;
\]
hence we have
\[
P\big(F_v; \phi(B\times \S^1)\big)=\int_{B\cap \{0<\theta_v<\pi/2\}}dr\int_{(\partial^*F_v)_r}\left(\frac{4r\theta_v'(r)}{p_{F_v}(r)}\right)^2d\mathcal{H}^0(x)=\int_B\sqrt{p_{F_v}^2(r)+(4r\theta'_v(r))^2}dr,
\]
where the last integral is over all $B$, since the integrand function vanishes almost everywhere on $\left\{\theta_v\in\left\{0,\pi/2\right\}\right\}$.\\
To conclude, it is then enough to show that 
\begin{equation*}
    \int_B\sqrt{p_{F_v}^2(r)+(4r\theta'_v(r))^2}dr\leq \int_B\sqrt{p_{E}^2(r)+(4r\theta'_v(r))^2}dr + 2 \mathcal{H}^1(\Gamma_E\cap B),
\end{equation*}
which is equivalent to
\begin{equation}\label{equiv strong ineq}
    \int_B\sqrt{p_{F_v}^2(r)+(4r\theta'_v(r))^2} - \sqrt{p_{E}^2(r)+(4r\theta'_v(r))^2}dr \leq 2 \mathcal{H}^1(\Gamma_E\cap B).
\end{equation}
Notice that, for all $r$ in $B\setminus \Gamma_E$, one has $p_{F_v}(r)\leq p_{E}(r)$ and so
\[
\int_{B\setminus \Gamma_E}\sqrt{p_{F_v}^2(r)+(4r\theta'_v(r))^2} - \sqrt{p_{E}^2(r)+(4r\theta'_v(r))^2}dr \leq 0.
\]
Hence, it is enough to prove that 
\[
\int_{B\cap \Gamma_E}\sqrt{p_{F_v}^2(r)+(4r\theta'_v(r))^2} - \sqrt{p_{E}^2(r)+(4r\theta'_v(r))^2}dr \leq 2 \mathcal{H}^1(\Gamma_E\cap B),
\]
which is implied by the fact that for all $r\in B\cap \Gamma_E$ we have
\begin{equation*}\label{rad-rad min 2}
    \sqrt{p_{F_v}^2(r)+(4r\theta'_v(r))^2} - \sqrt{p_{E}^2(r)+(4r\theta'_v(r))^2} \leq 2.
\end{equation*}
To obtain this last inequality, we simply notice that, if $r$ is in $B\cap \Gamma_E$, then $p_E(r)=2$, $p_{F_v}(r)=4$, and the inequality $\sqrt{16+t}-\sqrt{4-t}\leq 2$ is trivially true for all $t\geq 0$.
Summarizing, we showed that \eqref{equiv strong ineq} holds, and, as noticed above, this gives \eqref{strong ineq sph rearr}.
In particular, under the hypothesis of disconnection of non trivial circular slices of $E$, one has $\mathcal{H}^1(\Gamma_E\cap B)=0$, which implies \eqref{ineq sph rearr}, thus concluding the proof of the theorem.\hfill \qedsymbol
\begin{remark}
    It is straightforward to notice that, if we know that for all $r$ in a positive-measure subset of $B\subseteq (0,+\infty)$, the slice $E_r$ is made of three or more arcs of circle, then \eqref{ineq sph rearr} is a strict inequality.
\end{remark}
Another important remark is to be made. 
One might wonder whether the result proved in this paper is valid also in higher dimension.
The answer to this question is negative. 
Indeed, in the proof of Theorem \ref{thm sph rearr}, the key observation was that, for all $r>0$ 
\[
\mathcal{H}^0((\partial F_v)_r)\leq\mathcal{H}^0(\partial^* E_r),
\]  
due to the hypothesis on the slices of $E$.
In higher dimension there is not an analogous hypothesis which can guarantee such property of the rearrangement.
To understand this, it is sufficient to take, for example, a set $E$ made of a ball centered at the origin with two thin cylinders of the same length $l$ but with different section-radii $r,\alpha r$, for some $\alpha\neq 1$ fixed.
Then, the rearrangement of $E$ is a set made of a ball centered at the origin with two thin cylinders of length $l$ and section-radius $\tilde{r}$ such that $\tilde{r}^{N-1}=(r^{N-1}+(\alpha r)^{N-1})/2$ (see Figure \ref{fig counterex double sph high dim}).
Hence, we can estimate the perimeter of $E$ and the perimeter of its rearrangement as 
\[
\begin{aligned}
    P(E)&=\omega_N+(N-1)\omega_{N-1}l(r^{N-2}+(\alpha r)^{N-2})+O(r^{N-1})\\
    P(F_v)&=\omega_N+2(N-1)\omega_{N-1}l\tilde{r}^{N-2}+O(\tilde{r}^{N-1}).
\end{aligned}
\]
But we can rewrite $P(F_v)$ in terms of $r,\alpha r$ and, by concavity of the map $t\mapsto t^{\frac{N-2}{N-1}}$, we obtain
\[
P(F_v)=\omega_N+2(N-1)\omega_{N-1}l\left(\frac{(\alpha r)^{N-1}+r^{N-1}}{2}\right)^\frac{N-2}{N-1}> P(E),
\]
if $r\ll 1\ll l$.
Hence, we can deduce that in dimension higher than $2$ the double spherical cap symmetrisation does not decrease the perimeter in general, even under the assumption that non-trivial spherical slices are disconnected.



\begin{thebibliography}{99}
\bibitem{A} L. Ambrosio, ``Corso Introduttivo alla Teoria Geometrica della Misura ed alle Superfici Minime'', Appunti Corsi Tenuti Docenti Sc., Scuola Normale Superiore, Pisa, 1997.
\bibitem{AFP} L. Ambrosio, N. Fusco \& D. Pallara, ``Functions of Bounded Variation and Free Discontinuty Problems'', Oxford Math. Monogr., The Clarendon Press, Oxford University Press, New York, 2000.
\bibitem{BCH} C. Bianchini, G. Croce \& A. Henrot, {\it On the quantitative isoperimetric inequality in the plane with the barycentric distance}, Ann. Sc. Norm. Super. Pisa Cl. Sci.(5) {\bf 24} (2023), 2477-2500.
\bibitem{B} T. Bonnesen, {\it Les problèmes des isopérimètres et des isépiphanes}, Gautier-Villars, Paris, 1929.
\bibitem{C} S. Campi, {\it Isoperimetric deficit and convex plane sets of maximum translative discrepancy}, Geom. Dedicata {\bf 43} (1992), 71-81.
\bibitem{CPS} F. Cagnetti, M. Perugini \& D. Stoger, {\it Rigidity for perimeter inequality under spherical symmetrisation}, Calc. Var. Partial Differential Equations {\bf 59} (2020), Paper No 139, 53 pp.
\bibitem{CL} M. Cicalese \& G.P. Leonardi, {\it A selection principle for the sharp quantitative isoperimetric inequality}, Arch. Ration. Mech. Anal. {\bf 206} (2012), 617-643.
\bibitem{DG} E. De Giorgi, {\it Sulla proprietà isoperimetrica dell'ipersfera, nella classe degli insiemi aventi frontiera orientata di misura finita}, Atti Accad. Naz. Lincei Mem. Cl. Sci. Fis. Mat. Natur. Sez. la (8) {\bf 5} (1958), 33-44.
\bibitem{Fu1} B. Fuglede, {\it Stability in the isoperimetric problem}, Bull. London Math. Soc. {\bf 18} (1986), 599-605.
\bibitem{Fu2} B. Fuglede, {\it Stability in the isoperimetric problem for convex or nearly spherical domains in $\R^n$}, Trans. Amer. Math. Soc. {\bf 314} (1989), 619-638.
\bibitem{GP} C. Gambicchia \& A. Pratelli, {\it The sharp quantitative barycentric isoperimetric inequality for bounded sets}, to appear on Ann. Sc. Norm. Super. Pisa Cl. Sci.(5), 2024.
\bibitem{K} B. Kawhol, ``Rearrangemets and Convexity of Level Sets in PDE'', Lecture Notes in Math., Vol. 1150, Springer-Verlag, Berlin, 1985.
\bibitem{M} F. Maggi , ``Sets of Finite Perimeter and Geometric Variational Problems. An Introduction to Geometric Measure Theory'', Cambridge Stud. Adv. Math., Vol. 135, Cambridge University Press, Cambridge, 2012. 
\bibitem{P} G. P\'{o}lya, {\it Sur la symétrisation circulaire}, C.R. Acad. Sci. Paris {\bf 230} (1950), 25-27.
\end{thebibliography}
\end{document}